\documentclass{amsart}
\usepackage{amssymb}
\usepackage{amsfonts}
\usepackage{amsmath}
\usepackage[legalpaper,bookmarks=true,colorlinks=true,linkcolor=blue,citecolor=blue]{hyperref}
\usepackage{graphicx}%
\usepackage{fancyhdr}
\usepackage{color}
\usepackage[mathlines]{lineno}
\usepackage{lscape}
\usepackage{epsfig}
\usepackage{natbib}
\usepackage{geometry}
\usepackage{tgbonum}
\usepackage[]{listings}

\fontfamily{qcr}\selectfont

\theoremstyle{plain}

\newtheorem{definition}{Definition}

\numberwithin{equation}{section}

\newcommand{\Bin}{\bigskip \noindent}

\newcommand{\Ni}{\noindent}

\begin{document}
\Large
\title[Independence Indicator of records]{Independence of the indicator functions of record values for Multivariate independent data}

\begin{abstract}  We consider a sequence of random vectors on \(\mathbb{R}^d, \ d\geq 1\). We consider the record values based on the simultaneous strict inequality of the coordinates. The indicator record variable (irv) of the j-th observation is the function that assigns the value 1 (one) if that observation is a record value and the null value otherwise. Here, we give a detailed and a thorough proof  that the indicator functions are independent,  whenever the data are themselves independent, not necessarily \textit{iid}, in \(\mathbb{R}^d, \ d\geq 1\). We compare that proof  with available proofs in dimension one. Indeed, in seminal works on records, in particular in Ahsanullah(2024), Nevzorov(2001), Resnick (1987), Ahsanullah and Nevzorov (2015), etc., the independence of record indicator functions is usually validated based on logical reasoning, and so, is not rigorously proved. This allows us to undertake a detailed and thorough proof for independent data, not necessity \textit{iid}, in \(\mathbb{R}^d, \ d\geq 1\). The proof includes leads for not even independent data.\\.\\  

\noindent $^{\dag}$ Gane Samb Lo.\\
LERSTAD, Gaston Berger University, Saint-Louis, S\'en\'egal (main affiliation).\newline
LSTA, Pierre and Marie Curie University, Paris VI, France.\newline
AUST - African University of Sciences and Technology, Abuja, Nigeria\\
gane-samb.lo@edu.ugb.sn, gslo@aust.edu.ng, ganesamblo@ganesamblo.net\\
Permanent address : 1178 Evanston Dr NW T3P 0J9,Calgary, Alberta, Canada.\\

\noindent El hadji Babou\\
LERSTAD, Gaston Berger University, Saint-Louis, S\'en\'egal (main affiliation).\newline
Imhotep International Mathematical Center (imh-imc), https://imhotepsciences.org\\
Emails : babou.elhadji@ugb.edu.sn\\

\end{abstract}
\maketitle

\newpage
\bigskip \Ni \textbf{Presentation of the authors}.\\

\Ni \textbf{Gane Samb Lo}, Ph.D., is a retired full professor from Universit\'e Gaston Berger (2023), Saint-Louis, SENEGAL. He is the founder and The Probability and Statistics Chair holder at: Imhotep International Mathematical Center (imho-imc), https://imhotepsciences.org.\\

\Ni \textbf{El hadji Babou}, Ph.D., is preparing a Ph.D thesis at LERSTAD, Universit\'e Gaston Berger (2023), Saint-Louis, SENEGAL. under the supervision of the first author. He is also a junior researcher at: Imhotep International Mathematical Center (imho-imc), https://imhotepsciences.org.\\


\section{Introduction}

\Ni Let $X_1, \ X_2, \cdots $ be a sequence of real-valued random variables defined on some probability space 
$(\Omega, \mathbb{A}, \mathbb{P})$. Records times are integers $j=1$, and any $j\geq 2$ such that $M_{j-1}=\max(X_1, \cdots , X_{j-1})<X_j$ and $X_j$ is called a record value. For any $j\geq 1$, the function $\zeta_j$ takes the values one if $X_j$ is a record value and the null value $0$ otherwise.\\

\Ni It is well established that the sequence $(\zeta_j)_{j\geq 1}$ are independent if the $X_j$'s are independent and identically ditributed (\textit{iid}). Such a property is usually stated in the study of the sequence $L(n)$ the number if records in the $n\geq 1$ first observation. Indeed, we have

$$
L(n)= \sum_{i=1}^{n} \zeta_i.
$$

\Ni Consequently, the asymptotic law of $L(n)$ should be derived from the classical theory of partial sums of independent data. See
 \cite{loeve}, \cite{ips-mfpt-ang}, \cite{feller2}, etc.\\

\Ni Then the complete proof of the independence of the $\zeta_j$'s is important for the use of the central limit theorem of sums of independent random variables. In the dimension one, authors have usually explained more than proved. Here are some examples. In \cite{ahsan-tp}, on page 40, the author showed only that the $\zeta_j$'s are uncorrelated. In \cite{ahsan-nevz}, page 48, the proof of that independence is proposed as exercise thought this: \textit{Based on equalities (3.3.1), (3.3.2) and the analogous equalities for any set of indicators show that if the $X_j$'s are independent and have the same continuous distribution function $F$ then the $\zeta_j$'s  are independent random variables}. In \cite{nevzorov}, page 104, the independence is also stated from previous formulas that are not explicitly posed as independence conditions. Up to our knowledge, this is not treated in \cite{resnick}.\\

\Ni This motivated us to undertake the broader subject of the independence of the $\zeta_j$'s in higher dimensions for independent data, not necessarily with identical distribution. Records as not as simple as in dimension one. Let us give a simple example. A record time $j$  in dimension one for continuous data is defined as $X_j=max(X_1,\cdots,X_n)$. In dimension $d\geq 2$, $max(X_1,\cdots,X_n)$ does not always exists since the order is truly partial. However, despite the facts that formulations and formulas are more complex, the way we use the independence of that data to conclude the proof is the same.\\

\Ni We will give our proof using one the most general independence conditions as in \cite{ips-mfpt-ang}, Chapter 2 as follows.\\

\begin{definition} \label{def} Random variables $X_1, \cdots, X_n$, $n\geq 2$ random variables taking values in some measurable space 
$(F, \ \mathcal{F})$ and defined a probability space $(\Omega, \ \mathcal{A}, \ \ \mathbb{P})$ \textit{are mutually independent} if and only if for any non-negative and measurable real-valued functions : $h_j \ : E \rightarrow \mathbb{R}$, we have

\begin{equation}
\mathbb{E}\biggr(\prod_{j=1}^{n}\ h_j(X_j)\biggr)= \prod_{j=1}^{n} \mathbb{E}\left(h_j(X_j)\right). \label{F1}
\end{equation}
\end{definition}

\Ni In that context, since each $\zeta_j$ takes only two values, zero and one, we have to prove:  $\forall h_j:\left(\{0,1\}, \mathcal{P}(\{0,1\})\right)\to \mathbb{R}, \ j=2,\cdots,n$ (Here all these functions are measurable),

\begin{equation}
\mathbb{E}\left(\prod_{j=2}^{n}h_j(\zeta_j)\right)=\prod_{j=2}^{n}\mathbb{E}\left(h_j(\zeta_j)\right). \label{TBP1}
\end{equation}

\section{Proof}

\Ni Let $\mathcal{A}_k$ the set of ordered of subsets of $\mathcal{S}_n=\{2,\cdots,n\}$ of $k$ distinct elements. For $k\geq 3$, we represent an element $r$ of $\mathcal{A}_k$ as $r=(r_1,\cdots,r_k)$. For $k=1$, we write $r=(r_1)$ and for $k=2$, we write $r=(r_1,r_2)$. For $k=0$, we have  $\mathcal{A}_0=A_0$ and this means that $A_0$ is the subspace where $X_1$ the unique  record value.\\

\Ni So the set of all ordered sets in $\mathcal{S}_n$ of distinct elements is

$$
\mathcal{A}=\sum_{k=0}^{n-1}\mathcal{A}_k
$$ 

\Bin Now, let $r\in \mathcal{A}_k$, and $r^c=(j\in \mathcal{A}_k\setminus r)$. Let us set $A_0$ as the event : there is no record value after $X_1$, i.e. $A_0=(\zeta_2=0, \cdots, \zeta_n=0)$. For $k\geq 1$, for $r \in \mathcal{A}_k$, $A_r$ is the event: record values are 
$\{X_{r_1}, \cdots, X_{r_k}\}$. i.e.,   

$$
A_r=\left(\zeta_j=1, j\in r\right)\cap\left(\zeta_j=0, j\in r^c\right).
$$

\Bin We have to prove that: $\forall h_j:\left(\{0,1\}, \mathcal{P}(\{0,1\})\right)\to \mathbb{R}, \ j=2,\cdots,n$ (Here all these functions are measurable),

\begin{equation}
\mathbb{E}\left(\prod_{j=2}^{n}h_j(\zeta_j)\right)=\prod_{j=2}^{n}\mathbb{E}\left(h_j(\zeta_j)\right). \label{TBP}
\end{equation}

\Bin Let us denote the data as $X^{(j)}= (X_1^{(j)}, \cdots, X_d^{(j)})^T$. We already say that we use all-coordinate order, that is 
$x= (x_1, \cdots, x_d)^T \textit{ less than \ \ $(<)$ } y= (y_1, \cdots, y_d)^T$ if and only if

$$
\forall i \in \overline{1, \ d}, \ \ x_i<y_i.
$$

\Bin For $r=\emptyset$, we have on  $A_0$:

\begin{eqnarray*}
A_0=B_1(X^{(1)})&=&\left\{\forall j \in \overline{2, \ n}, \ \neg \left(X^{(j)}>X^{(1)}\right) \right\}\\
&=&\left\{\forall j \in \overline{2, \ n}, \left(\exists \ i \in \overline{1, \ d}\right), \ \left(X_i^{(j)}\leq X_i^{(1)}\right) \right\}
\end{eqnarray*}

\Bin Let us define for $r \in \mathcal{A}_k$, for $2\leq \ell \leq k$,

\begin{eqnarray*}
A(X^{(r_{\ell-1})}, X^{(r_\ell)})&=&\left\{\forall j \in \overline{r_{\ell-1}+1, r_{\ell}-1}, \neg \left(X^{(j)}>X^{(r_{\ell-1})}\right)  \textit{ and } 
 \left(X^{(r_{\ell-1})}<X^{(r_\ell)}\right) \right\}
\end{eqnarray*}

\Bin and

\begin{eqnarray*}
B_k(X^{(r_k)})&=&\left\{\forall j \in \overline{r_{k}+1, n}, \ \ \neg \left(X^{(j)}>X^{(r_k)}\right) \right\}
\end{eqnarray*}

\noindent with the convention that $B_k(X^{(r_k)})=\Omega$ if $\overline{r_{k}+1, n}=\emptyset$, i.e., $r_k=n$.\\

\noindent Then for let us define for $r \in \mathcal{A}_k$, we denote $X^{[r]}=(X^{r_1}, \cdots, X^{r_k})^T$ and as well 
$x^{[r]}=(x^{r_1}, \cdots, x^{r_k})^T$ from ()
 
\begin{eqnarray*}
A_r&=&A_r\left(X^{[r]}\right)\\
&=&B_1\left(X^{(1)}\right)\cap A\left(X^{(1)}, X^{(2)}\right) \cap A\left(X^{(2)}, X^{(3)}\right) \cap \ \cdots \ \cap A\left(X^{(r_{k-1}}), X^{(r_{k})}\right)  \cap B_k\left(X^{(r_k)}\right).
\end{eqnarray*}

\Bin Now, for 

$$
Z=\left(\prod_{j=2}^{n}h_j(\zeta_j)\right)
$$

\Bin we have

\begin{equation*}
\mathbb{E}Z= \sum_{k=0}^{n-1} \sum_{A_r \in \mathcal{A}_k} \mathbb{E}\left(Z/A_r\right) \mathbb{P}(A_r).
\end{equation*}

\Bin On $A_r$, we have

$$
Z=\prod_{j \in r} h_j(1) \times \prod_{j \in r^c} h_j(0).
$$

\Bin We notice that :\\

\Ni \textbf{Remark (R1)} \label{R1}, for $x \in \mathbb{R}^d$, the events  $B_1(x^{(1)})$, $B_k(x^{r_k})$ and the $A(x^{r_{\ell-1}}, x^{r_{\ell}})$ are independent since these random variables are functions of distinct blocks of the independent data. Then, with the re-writing of indicator function as $1_{\circ}=I({\circ})$, we have

\begin{eqnarray}
 \ \ \mathbb{E}I(A_r)&=&\int \mathbb{E}\left(I(A_r)\biggr/\left(X^{[r]}=x^{[r]}\right)\right) \ d\mathbb{P}\left(x^{[r]}\right)\notag\\
\ \ &=&\int \mathbb{E}\left(I(B_1(X^{(1)})) \times    \prod_{\ell=2}^{k} I(A(X^{r_{\ell-1}}, X^{r_{\ell}})) \times I(B_k(X^{r_k}))\biggr
/\left(X^{[r]}=x^{[r]}\right)\right) \ d\mathbb{P}\left(x^{[r]}\right) \notag\\
\ \ &=&\int \mathbb{E}\left(I(B_1(x^{(1)})) \times    \prod_{\ell=2}^{k} I(A(x^{r_{\ell-1}}, x^{r_{\ell}})) \times I(B_k(x^{r_k}))\right)
 \ d\mathbb{P}\left(x^{[r]}\right) \notag\\
\ \ &=&\int \mathbb{E}\left(I(B_1(x^{(1)}))\right) \ d\mathbb{P}\left(x^{[r_k]}\right) \times  \int  \prod_{\ell=2}^{k} \mathbb{E}\left(I(A(x^{r_{\ell-1}}, x^{r_{\ell}}))\right) \ d\mathbb{P}\left(x^{[r]}\right)\\
&\times& \int \mathbb{E}\left(I(B_k(x^{r_k}))\right)  \ d\mathbb{P}\left(x^{[r_k]}\right) \notag\\
&=&\mathbb{E}\left( I(B_1(X^{(1)}))\right) \times    \prod_{\ell=2}^{k} \mathbb{E}\left(I(A(X^{r_{\ell-1}}, X^{r_{\ell}}))\right) \times \mathbb{E}\left(I(B_k(X^{r_k}))\right).\label{L5}
\end{eqnarray}

\Bin Notice that in Line 5, we used the independence described just above and in Line 4, we have an integration by a product measures $d\mathbb{P}\left(x^{[r]}\right)= d\mathbb{P}\left(x^{r_1}\right) \cdots d\mathbb{P}\left(x^{r_k}\right)$ associated the convolution product of the $X^{r_j}$s.\\

\Ni Now, still on $A_r$, we have respectively on $B_1(X^{(1)})$, on the $A(X^{r_{\ell-1}}, X^{r_{\ell}})$'s and on $B_k(X^{r_k})$, the products

\begin{eqnarray}
&&v_0(\zeta)=\prod_{j=2}^{n}  h_j(\zeta_j)=\prod_{j=2}^{n}  h_j(0)\\
&&v_\ell(\zeta)=\prod_{j=r_{\ell-1}+1}^{r_{\ell}}  h_j(\zeta_j)=\biggr(\prod_{j=r_{\ell-1}+1}^{r_{\ell}-1} \   h_j(0)\biggr) \ h_{r_\ell}(1)\\
&&v_{k+1}(\zeta)=\prod_{j=r_{k+1}}^{n}  h_j(\zeta_j)=\prod_{j=r_{k+1}}^{n}  h_j(0).
\end{eqnarray}

\Bin So, we can see that $\mathbb{E}\prod_{j=2}^{n} \ h(\zeta_k)$ takes the values $v_{\circ}$ on the elements of decomposition of $\Omega$, i.e., on $A_0$, $A_r \in \mathcal{A}_k, \ \overline{k\in 1,n-1}$. So it is a piece-wise constant function since we have

\begin{equation}
\mathbb{E}\prod_{j=2}^{n} \ h(\zeta_k)=v_0(\zeta) \mathbb{E}(I(A_0))+  \sum_{k=1}^{n-1}\sum_{A_r \in \mathcal{A}_k} v_0(\zeta)\left(\prod_{\ell=2}^{k}v_\ell(\zeta)\right) v_{k+1}(\zeta)\mathbb{E}(I(A_r)) . \label{L10}
\end{equation}

\Bin Now let us combine \eqref{L5} and \eqref{L10} to conclude. Indeed, we have

\begin{eqnarray*}
\mathbb{E}Z &=& \mathbb{E}( Z \ I(A_0)\\
						&+& \sum_{k=1}^{n-1} \sum_{A_r \in \mathcal{A}_k} \mathbb{E}\left(Z \ I({A_r})\right)\\
						&=:&\Sigma_1 + \Sigma_2.
\end{eqnarray*}

\Bin Let us first focus on 

\begin{eqnarray*}
\Sigma_2&=& \sum_{k=1}^{n-1} \sum_{A_r \in \mathcal{A}_k} \mathbb{E}\biggr(Z \ \times \left(I(B_1(X^{(1)}))\right) \times  \left( \prod_{\ell=2}^{k} I(A(X^{r_{\ell-1}}, X^{r_{\ell}}))\right) \times \left(I(B_k(X^{r_k})\right)\biggr)\\
&=& \sum_{k=1}^{n-1} \sum_{A_r \in \mathcal{A}_k} \int \mathbb{E}\biggr(Z \ \times \biggr(I(B_1(X^{(1)}))\biggr) \times  \biggr( \prod_{\ell=2}^{k} I(A(X^{r_{\ell-1}}, X^{r_{\ell}}))\biggr) \\
&\times& \biggr(I(B_k(X^{r_k})\biggr) \biggr/ (X^{[r]}=x^{[r]}) \biggr) \ d\mathbb{P}_{X^{[r]}}(x^{[r]})\\
&=& \sum_{k=1}^{n-1} \sum_{A_r \in \mathcal{A}_k} \int \mathbb{E}\biggr(Z \ \times \biggr\{ I(B_1(x^{(1)}))\biggr\} \times  \biggr( \prod_{\ell=2}^{k} I(A(x^{r_{\ell-1}}, x^{r_{\ell}}))\biggr) \\
&\times& \biggr\{I(B_k(x^{r_k})) \biggr\} \ d\mathbb{P}_{X^{[r]}}(x^{[r]})\\
&=& \sum_{k=1}^{n-1} \sum_{A_r \in \mathcal{A}_k} \int \mathbb{E}\biggr\{\biggr(v_0(\zeta) \ I(B_1(x^{(1)}))\biggr) \times  
\biggr(\prod_{\ell=2}^{k} \ \biggr(v_{\ell} \ I(A(x^{r_{\ell-1}}, x^{r_{\ell}}))\biggr) \\
&\times& \biggr(v_{k+1}(\zeta) \ I(B_k(x^{r_k})) \biggr)\biggr\} \ d\mathbb{P}_{X^{[r]}}(x^{[r]})\\
\end{eqnarray*}

\Bin by the independence described in Remark \ref{R1} and the constancy of the $v_{\circ}$s, 

\newpage
\begin{eqnarray}
\Sigma_2&=& \sum_{k=1}^{n-1} \sum_{A_r \in \mathcal{A}_k} \int  \mathbb{E}\biggr(v_0(\zeta) \ I(B_1(x^{(1)})) \ d\mathbb{P}_{X^{[r]}}(x^{[r]})\biggr) \times  
\biggr(\prod_{\ell=2}^{k} \ \mathbb{E}\biggr(v_{\ell}(\zeta) \ I(A(x^{r_{\ell-1}}, x^{r_{\ell}})) \ d\mathbb{P}_{X^{[r]}}(x^{[r]})\biggr)\notag \\
&\times& \mathbb{E}\biggr(v_{k+1}(\zeta) \ I(B_k(x^{r_k})) \ d\mathbb{P}_{X^{[r]}}(x^{[r]})\biggr) \notag\\
&=& \sum_{k=1}^{n-1} \sum_{A_r \in \mathcal{A}_k} \int \mathbb{E}\biggr(v_0(\zeta)\ I(B_1(x^{(1)})) \ d\mathbb{P}_{X^{[r]}}(x^{[r]})\biggr) \times  
\biggr(\prod_{\ell=2}^{k} \ \mathbb{E}\biggr( v_{\ell}(\zeta) I(A(x^{r_{\ell-1}}, x^{r_{\ell}})) \ d\mathbb{P}_{X^{[r]}}(x^{[r]})\biggr) \notag\\
&\times& \mathbb{E}\biggr(v_{k+1}(\zeta) \ I(B_k(x^{r_k})) \ d\mathbb{P}_{X^{[r]}}(x^{[r]})\biggr) \notag\\
&=& \sum_{k=1}^{n-1} \sum_{A_r \in \mathcal{A}_k} \prod_{j=2}^{n}\biggr\{\mathbb{E}h(\zeta_k)\biggr\} \  \int \mathbb{E}\biggr(\ I(B_1(x^{(1)})) \ d\mathbb{P}_{X^{[r]}}(x^{[r]})\biggr) \label{II1}\\
&\times&  
\int\biggr(\prod_{\ell=2}^{k} \ \mathbb{E}\biggr( \ I(A(x_{r_{\ell-1}}, x_{r_{\ell}})) \ d\mathbb{P}_{X^{[r]}}(x^{[r]})\biggr) 
\times\int\mathbb{E}\biggr( \ I(B_k(x_{r_k})) \ d\mathbb{P}_{X^{[r]}}(x^{[r]})\biggr) \label{l1}\\
&=& \prod_{j=2}^{n}\biggr\{\mathbb{E}h(\zeta_k)\biggr\} \sum_{k=1}^{n-1} \sum_{A_r \in \mathcal{A}_k} \
\mathbb{E}\biggr(\ I(B_1(X^{(1)})) \prod_{\ell=2}^{k} \ \mathbb{E}\biggr( \ I(A(X^{r_{\ell-1}}, X^{r_{\ell}})) I(B_k(X^{r_k}))\biggr) \label{l2} \\
&=& \prod_{j=2}^{n}\biggr\{\mathbb{E}h(\zeta_k)\biggr\} \sum_{k=1}^{n-1} \sum_{A_r \in \mathcal{A}_k} \mathbb{P}(A_r(X)) \label{l3}
\end{eqnarray}

\Bin Notice first that in Line \ref{II1}, the expectation%
$$
\prod_{j=2}^{n}\biggr\{\mathbb{E}h(\zeta_k)\biggr\}
$$

\Ni is, actually, a constant on each $A_r$, and so it is piece-wise contant function on $\Omega$. Secondly, we are able to regroup the expectations after the double summation ni one because of the independence and next use the continuous Bayes formula to get the final unconditional expectation.\\

\Ni Next, on $\Sigma_1$, all the $h_j(\zeta)$, $j\geq 2$ are $h_j(0)$, and then

\begin{eqnarray}
\Sigma_1= \prod_{j=2}^{n}\biggr\{\mathbb{E}h(\zeta_k)\biggr\} \mathbb{P}(A_0(X)). \label{l4}
\end{eqnarray}

\noindent By putting \eqref{l3} and \eqref{l4} together we conclude that \eqref{TBP} is well established. The proof is complete.\\

\label{fin-art}
\end{document}